\documentclass[12pt]{amsart}
\newcommand{\ZZ}{\ensuremath{\mathbb{Z}}}

\usepackage{amssymb,a4wide,epsfig,psfrag}
\usepackage{url}
\newtheorem{theorem}{Theorem}
\newtheorem{proposition}[theorem]{Proposition}

\newfont{\Ar}{msam10}

\newcommand{\cK}{\ensuremath{{\mathcal K}}}
\DeclareMathOperator{\Stab}{Stab}
\DeclareMathOperator{\Cay}{Cay}
\DeclareMathOperator{\Syl}{Syl}
\newcommand{\RR}{\ensuremath{\mathbb{R}}}

\title{Wythoff polytopes and low-dimensional homology of Mathieu groups}

\author{Mathieu Dutour Sikiri\'c}
\address{M. Dutour Sikiri\'c, Rudjer Boskovi\'c Institute, Bijenicka 54, 10000 Zagreb, Croatia}
\email{mdsikir@irb.hr}

\author{Graham Ellis}
\address{G. Ellis, Mathematics Department, National University of Ireland, Galway}
\email{graham.ellis@nuigalway.ie}

\thanks{First author has been supported by Marie Curie fellowship MTKD-CT-2006-042685 and by the Croatian Ministry of Science, Education and Sport under contract 098-0982705-2707.}

\begin{document}
\maketitle

\begin{abstract}
We describe two methods for computing 
 the low-dimensional 
integral homology of the Mathieu simple groups and use them to make computations such as 
$H_5(M_{23},\ZZ)=\ZZ_7$ and $H_3(M_{24},\ZZ)=\ZZ_{12}$. 
One method works via Sylow subgroups.
The other method uses a Wythoff polytope and perturbation techniques to produce an explicit free $\ZZ M_n$-resolution.
Both methods apply in principle to arbitrary finite groups.
\end{abstract}

\section{Introduction}
We describe two methods for computing the integral homology
for the Mathieu simple groups presented on Table \ref{TableSmallHomology}.
The first homology $H_1(G,\ZZ)$ is 
trivial for any simple group and so is omitted from the table
(see \cite{brown} for an exposition of relevant facts on group homology). 
The second homology of Mathieu groups is well-known \cite{mazet}.
A computer method for the second homology of
a permutation group was illustrated on the Mathieu groups
$M_{21}$ and $M_{22}$ in \cite{holt}.
The mod $p$ cohomology 
$H^*(G,\mathbb F_p)$ is now known for all Mathieu groups 
except $M_{24}$ \cite{webb,ademmaginnismilgram,ademmilgram,milgram}.
With the help of the Bockstein spectral sequence it is, in principle,
 possible to obtain
integral homology from mod $p$ 
cohomology ($p$ ranging over the prime divisors of the  group order),
 though the details can be difficult.
For example, the calculation of $H_n(M_{23},\ZZ)$ was obtained
in this way for $1\le n\le 6$ by Milgram \cite{milgram} and 
provided the first example of a non-trivial finite group with 
trivial integral homology in dimensions $\le 3$.
It seems that the mod $p$ cohomology of $M_{24}$ is not
known for all primes $p$ (see \cite{green} for the case $p=3$)
and so we can assign the status of a new theorem to the following result.

\begin{theorem}\label{h3m24}
$H_3(M_{24},\ZZ) =\ZZ_{12}$ and $H_4(M_{24},\ZZ) = 0$.
\end{theorem}
This result (and other table entries) can be obtained from
the
{\sc hap} homological algebra 
package \cite{hap} for the {\sc gap} computational algebra system \cite{gap}
using (variants of) the following command.

\begin{table}
\begin{equation*}
\begin{array}{l|l|l|l|l}
  G &H_2(G,\ZZ) &H_3(G,\ZZ) &H_4(G,\ZZ) &H_5(G,\ZZ)\\
\hline
M_{11} &0 &\ZZ_8 &0 &\ZZ_2\\
M_{12} &\ZZ_2 &\ZZ_6\oplus \ZZ_8 &\ZZ_3 &(\ZZ_2)^3\\
M_{21} &\ZZ_4 \oplus \ZZ_{12} & \ZZ_{5} &0 &(\ZZ_2)^4\oplus \ZZ_4 \oplus \ZZ_7\\
M_{22} &\ZZ_{12} & 0 &0 &  \ZZ_2 \oplus \ZZ_2 \oplus   \ZZ_7\\
M_{23} & 0 & 0 &0 &\ZZ_7\\
M_{24} &0 &\ZZ_{12} & 0 &  (\ZZ_2)^a 
\oplus (\ZZ_4)^b \oplus \ZZ_7
\end{array}
\end{equation*}
\caption{Low dimensional homology of Mathieu groups with $0\le a\le 53$ and $0\le b\le 1$.}
\label{TableSmallHomology}
\end{table}

\begin{verbatim}
gap> GroupHomology(MathieuGroup(24),3);
gap> [ 4, 3 ]
\end{verbatim}

\noindent The algorithm underlying this command is explained
in Section \ref{alg}.
The current implementation is unable to determine the integers
$a,b$ in Table \ref{TableSmallHomology}
though it does establish the ranges
$0\le a\le 53$, $0\le b\le 1$. 

Abelian invariants of a (co)homology group are
 the easiest cohomological information to access. More 
difficult information would be, for example, explicit cocycles $G^n\rightarrow A$ corresponding to
cohomology classes in $H^n(G,A)$. 
Explicit cocycles are constructed in 
{\sc hap} using the induced chain map $B_{\ast}^G \rightarrow R_{\ast}^G$
from the bar resolution $B_\ast^G$
to an explicit small free $\ZZ G$-resolution
$R_{\ast}^G$ of $\ZZ$.
In Sections 3-5 we explain how the Wythoff polytope construction
can be used to produce such a resolution $R_{\ast}^G$.
This resolution provides an alternative computation of $H_3(M_{24},\ZZ)$.

In Section 6 we determine 
the $p$-part $H_n(M_{m},\ZZ)_{(p)}$ of the integral homology of
the Mathieu groups for $n\ge 1$ and primes $p\ge 5$.
For $p\in \{5,7,11,23\}$ the $p$-part is either trivial or $\ZZ_p$;
it is trivial for all other primes $p\ge 5$.
Table \ref{ZPvaluesNonZeroHomology} lists the values of $n$ for which the
$p$-part is non-trivial.

Although the paper focuses on Mathieu groups, the techniques are
applicable in principle to arbitrary finite groups.
In some cases the Wythoff polytopal method is a significantly faster
method for computing the homology groups.

\section{Algorithm underlying the {\sc hap} function}\label{alg}
Given a group $G$, a {\em free $\ZZ G$-resolution}  of the trivial module
$\ZZ$ is an exact sequence
\begin{equation*}
0\leftarrow \ZZ \leftarrow R_0^G\leftarrow R_1^G \leftarrow \dots \leftarrow 
R_k^G \leftarrow \dots
\end{equation*}
of free $\ZZ G$-modules  $R_i^G$.
A previous paper \cite{jsc} describes an algorithm for computing
free $\ZZ G$-resolutions  for finite $G$.
This has now been implemented as part of the {\sc hap} package.
It takes as input a finite group $G$ and a positive integer $n$.
It returns:
\begin{itemize}
\item The rank of the $k$th module $R_k^G$ in a
free $\ZZ G$-resolution $R_\ast^G$ ($0\leq k\le n$).
\item The image of the $i$th free $\ZZ G$-generator of $R_k^G$ under the 
boundary homomorphism $d_k\colon R_k^G\rightarrow R_{k-1}^G$ ($1\leq k\le n$).
\item The image of the $i$th free $\ZZ$-generator of $R_k^G$ under a  
contracting homotopy $h_k\colon R_k^G\rightarrow R_{k+1}^G$ ($0\leq k\le n-1$).
\end{itemize} 
The contracting homotopies $h_k$ satisfy, by definition, $h_kd_{k+1}+d_{k+2}h_{k+1}=1$
and need to be specified on a set of free 
Abelian group generators of $R_k$ since they are not $G$-equivariant.
The homotopy can be used to make constructive the following 
frequent element of choice.
$$\begin{minipage}{4.7in}\it{
For  $x\in \ker(d_k\colon R_k^G\rightarrow R_{k-1}^G)$ choose an element $\tilde x \in R_{k+1}^G$ such that $d_{k+1}(\tilde x) = x$.}
\end{minipage}$$
One sets $\tilde x = h_k(x)$.
In particular, for any group homomorphism $\phi\colon G\rightarrow G'$,
the homotopy allows one to define an induced 
$\phi$-equivariant chain map $\phi_\ast\colon R_\ast^G\rightarrow R_\ast^{G'}$. 

\medskip
The algorithm in \cite{jsc} can only handle fairly small groups. For example, the {\sc hap} implementation
takes $20$ seconds on a  2.66GHz Intel PC with 2G of memory 
to compute eight terms of a free $\ZZ G$-resolution $R_\ast^G$
for the symmetric group $G=S_5$; the $\ZZ G$-rank of $R_8^G$ is $115$.
However, for any group $G$ there is a surjection
$$H_n(\Syl_p,\ZZ) \rightarrow H_n(G,\ZZ)_{(p)}$$ from the homology of 
a Sylow $p$-subgroup $\Syl_p=\Syl_p(G)$
onto the $p$-part of the homology of $G$.
For a Sylow $p$-subgroup $P$ there is a  description 
of the kernel of the surjection $H_n(P,Z) \rightarrow H_n(G,Z)_{(p)}$
due to Cartan and Eilenberg \cite{ce}.
It is generated by elements
\begin{equation*}
\phi_K(a)  -  \phi_{xKx^{-1}}(a)
\end{equation*}
where $x$ ranges over the double coset representatives of $P$ in $G$,
$K=P \cap xPx^{-1}$, the homomorphisms $\phi_K$, 
$\phi_{x^{-1}Kx}\colon H_n(K,\ZZ) \rightarrow H_n(P,\ZZ)$ are induced by the inclusion $K\rightarrow P, k\mapsto k$ and the conjugated inclusion
$K\rightarrow P, k\mapsto x^{-1}kx$, 
and $a$ ranges over the generators of  $H_n(K,\ZZ)$.
Thus, the homology of a large finite group $G$ can be computed 
from free resolutions (with specified contracting homotopy) for each
of its Sylow subgroups.
Our implementation of the algorithm in \cite{jsc} can be used
to produce six terms of free $\ZZ(\Syl_p)$-resolutions for all
Sylow subgroups $\Syl_p$ of all Mathieu groups except $M_{24}$.
The Sylow subgroup $\Syl_2(M_{24})$ has order $1024$ and  requires 
a specific application of a general technique. 

\begin{table}
\begin{equation*}
\begin{array}{l|l|l|l|l|l|l}
H_n(M_m,\ZZ)_{(p)}=\ZZ_p &m=11 &m=12 &m=21 &m=22 &m=23 &m=24\\\hline
                          &n=   &n=   &n=   &n=   &n=   &n=\\\hline
p=5           &8k-1   &8k-1   &4k-1  &8k-1   &8k-1   &8k-1\\
p=7           &-      &-      &6k-1  &6k-1   &6k-1   &6k-1\\
p=11          &10k-1  &10k-1  &-     &10k-1  &10k-1  &20k-1\\
p=23          &-      &-      &-     &-      &22k-1  &22k-1\\
\end{array}
\end{equation*}
\caption{Values of $n$ expressed in term of $k\geq 1$ such that $H_n(M_m,\ZZ)_{(p)}=\ZZ_p$}
\label{ZPvaluesNonZeroHomology}
\end{table}

\medskip
To explain the technique suppose that $G$ is a group, possibly infinite,
for which we have some $\ZZ G$-resolution of $\ZZ$
\[C_\ast \colon \cdots\rightarrow C_n\rightarrow C_{n-1} \rightarrow \cdots \rightarrow C_0 \rightarrow \ZZ .\]
but that $C_\ast$ is not free.
Suppose that for each $m$ we have a free $\ZZ G$-resolution of
the module $C_m$
\[D_{m*} \colon \rightarrow D_{m,n} \rightarrow D_{m,n-1} \rightarrow
\cdots \rightarrow D_{m,0} \rightarrow C_m .\]

\begin{theorem}\label{wall}\cite{wall}
There is a free $\ZZ G$-resolution $R_{\ast}^G\rightarrow \ZZ$ with
\[ R_n^G =\bigoplus_{p+q=n} D_{p,q}\]
\end{theorem}
The proof of this theorem of C.T.C. Wall can be made constructive 
by using contracting homotopies on the resolutions $D_{m*}$. 
Furthermore, a contracting homotopy on $R_\ast^G$ can be constructed by a 
formula involving contracting homotopies on the $D_{m*}$ and on $C_\ast$.
Details are given in \cite{ellisharrisskoldberg}.

\medskip
Suppose now that $N$ is a normal subgroup of $G$ and that 
 $C_\ast$ is a free $\ZZ(G/N)$-resolution. Then, regarding $C_\ast$
as a $\ZZ G$-resolution, 
each free $\ZZ G$-generator of $C_m$ is stabilized 
by $N$. Any free $\ZZ N$-resolution of $\ZZ$ can be
 used to construct a free $\ZZ G$-resolution $D_{m*}$ of
$C_m$. Thus, 
using Theorem \ref{wall}, we can construct a free $\ZZ G$-resolution 
$R_\ast^G$ from a free
$\ZZ N$-resolution $R_\ast^N$ and free
$\ZZ(G/N)$-resolution $R_\ast^{G/N}$. 
The constructed resolution is often referred to as a {\it twisted tensor product} and denoted by $R_\ast^G=R_\ast^N\tilde\otimes R_\ast^{G/N}$.
 
\medskip
This twisted tensor product has been implemented in {\sc hap} and can be
used to provide free resolutions for the Sylow subgroup $\Syl_p(M_{24})$.
Since $|M_{24}|= 2^{10}\cdot 3^3\cdot 5\cdot 7\cdot 11\cdot 23$
the non-cyclic Sylow subgroups occur only for $p=2$, $3$.
Their low-dimensional integral homology can be computed
using {\sc hap} and is given in Table \ref{HomologyM24_p23}.

\begin{table}
\begin{equation*}
\begin{array}{l|l|l|l|l|l}
p & H_1(\Syl_p,\ZZ) & H_2(\Syl_p,\ZZ) & H_3(\Syl_p,\ZZ) & H_4(\Syl_p,\ZZ) & H_5(\Syl_p,\ZZ)\\\hline
&&&&&\\[-3mm]
2 &(\ZZ_2)^4 &(\ZZ_2)^8 &(\ZZ_2)^{11} \oplus (\ZZ_4)^{6} &(\ZZ_2)^{32} &(\ZZ_2)^{52} \oplus \ZZ_4\\
3 &(\ZZ_3)^2 &(\ZZ_3)^2 &(\ZZ_3)^4  &(\ZZ_3)^3 &(\ZZ_3)^4\oplus \ZZ_9
\end{array}
\end{equation*}
\caption{Low dimensional homology of Sylow subgroups of $M_{24}$ for $p=2$, $3$}
\label{HomologyM24_p23}
\end{table}

In degrees $n=5$ the current version of {\sc hap} 
fails to determine the image of $H_n(\Syl_2,\ZZ)$ in $H_n(M_{24},\ZZ)$.
It succeeds in constructing the image as a finitely presented group
but fails to determine the group from this presentation.
This failure should be resolved in a future release of {\sc hap}.

\medskip
The remainder of the paper is aimed at constructing small
free resolutions for large groups such as $M_{24}$.

\section{Orbit polytopes}

Suppose that a finite group $G$ acts linearly on $\RR^n$. For a vector
$v\in \RR^n$ we consider the convex hull
$$P=P(G,v) = {\rm Conv}(v^g : g\in G)$$
of the orbit of $v$ under the action of $G$.
The polytope $P$ has a natural cell structure with respect
to which we can consider the cellular chain complex $C_{\ast}(P)$.
The action of $G$ on $\RR^n$ induces an action of $G$ on
$C_\ast(P)$ and we can view $C_\ast(P)$ as a chain complex of $\ZZ G$-modules.  
Since $P$ is contractible we have 
$H_i(C_\ast(P)) =0$ for all $i\ge 1$ and $H_0(C_\ast(P)) =\ZZ$.
Furthermore, if the polytope is of dimension $m$ then 
$H_0(C_\ast(P)) \cong \ZZ \cong C_m(P)$.
So there is a homomorphism $C_0(P) \rightarrow C_{m-1}(P)$
which can be used to splice together infinitely many copies
of $C_\ast(P)$ to form an infinite $\ZZ G$-resolution
$$\cdots \rightarrow C_1 \rightarrow C_0 \rightarrow C_{m-1} \rightarrow \cdots \rightarrow C_2 \rightarrow  C_1 \rightarrow C_0 \rightarrow \mathbb Z $$
of the trivial $\ZZ G$-module $\ZZ$.
In principle one can use Theorem \ref{wall} to convert $C_\ast$
to a free $\ZZ G$-resolution.
Precise details are given in \cite{ellisharrisskoldberg}.
To put this idea into practice one requires:
\begin{enumerate}
\item The face lattice of the orbit polytope $P(G,v)$.
\item For each orbit of cell $e$ in $P(G,v)$, the subgroup $\Stab(G,e) \le G$ of elements that stabilize $e$ globally.
\item A free $\ZZ \Stab(G, e)$-resolution $R_\ast^{\Stab(G, e)}$ for each stabilizer $\Stab(G,e)$.
\end{enumerate} 
Assuming that the stabilizer groups $\Stab(G,e)$ are reasonably small,
resolutions $R_\ast^{\Stab(G,e)}$ are readily obtained from {\sc hap}'s
implementation of the algorithm in \cite{jsc}.
Thus, to convert $C_\ast$ to a free $\ZZ G$-resolution, we must focus on requirements (1) and (2). 

\bigskip
  One could  use computational geometry software
 such as Polymake \cite{polymake} to determine the
 combinatorial structure of $P(G,v)$ for small groups $G$. 
 For instance, any permutation group  $G\le S_n$ acts 
on $\RR^n$ by
$\pi(x_1,\ldots,x_n) = (x_{\pi^{-1}(1)},
 \ldots, x_{\pi^{-1}(n)})$ for $\pi\in G$. 
In particular, the Mathieu group $M_{10}$ of order $720$, generated by 
$\pi_1=( 1, 9, 6, 7, 5)( 2,10, 3, 8, 4)$ and
$\pi_2=( 1,10, 7, 8)( 2, 9, 4, 6)$, acts on $\RR^{10}$.
For the vector  
$v=(1,2,3,4,5,6,7,8,9,10)$ the polytope
$P(M_{10},v)$ is $9$-dimensional with $720$ vertices each of
degree $632$.
The polytope thus has $227520$ edges.

\section{Orbit polytopes of finite reflection groups}\label{Orbit}

Let $W$ be a finite reflection group generated by
a simple system of Euclidean reflections
$S=\{s_1, \ldots, s_n\}$.
For each reflection $s\in S$ let
$H_s$ denote the corresponding reflecting hyperplane
and $\Delta$ the fundamental simplex for $S$.
The Coxeter-Dynkin reduced diagram is the graph on $S$ with two reflections
adjacent if they do not commute.
Fix a subset $\emptyset\subsetneq V\subseteq S$.
The {\em type} $T=t(v) \subset S$ of a point $v\in \Delta$ is the set
of $s\in S$ such that $v\notin H_{s}$.
Choose a point $v$ of type $V$.
Let $P(W;V,v)$ denote the $n$-dimensional polytope formed by
the convex hull of the orbit of $v$ under the action of $W$.

\begin{figure}
\begin{center}
\begin{minipage}{5.2cm}
\centering
\resizebox{34mm}{!}{\rotatebox{0}{\includegraphics[bb=172 48 442 351,clip]{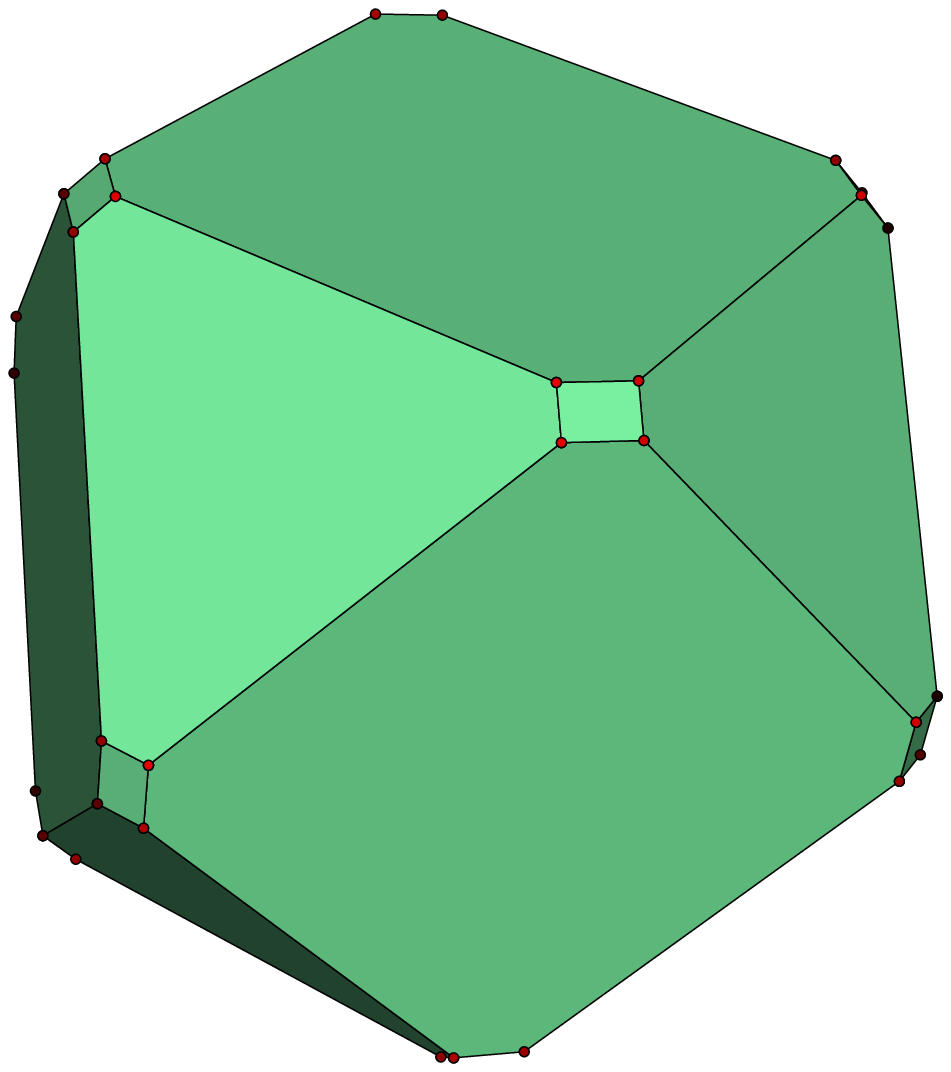}}}\par
a) The polytope $P(\mathsf{B}_3;\{s_1, s_2, s_3\}, v)$
\end{minipage}
\begin{minipage}{5.2cm}
\centering
\resizebox{34mm}{!}{\rotatebox{0}{\includegraphics[bb=172 40 442 354,clip]{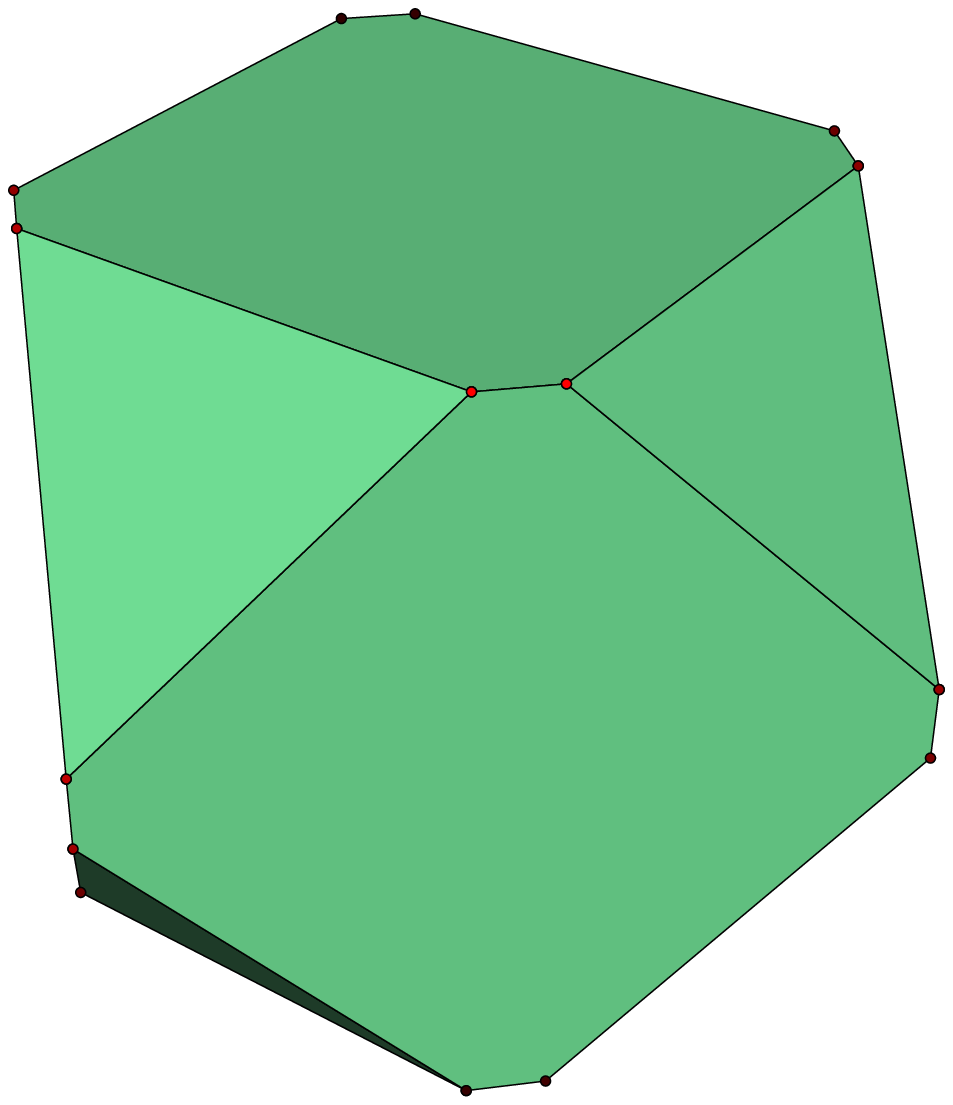}}}\par
b) The polytope $P(\mathsf{B}_3;\{s_2, s_3\}, v)$
\end{minipage}
\end{center}
\caption{Two Wythoff polytopes constructed from $D_3$}
\label{ExampleWythoff}
\end{figure}

As an example, consider the $3$-dimensional 
reflection group $W=\mathsf{B}_3$ generated by reflections 
$s_1,s_2,s_3$ where $(s_1s_2)^3=1$, $(s_1s_3)^2=1$ and $(s_2s_3)^4=1$. 
For $V=\{s_1, s_2, s_3\}$ and vector
$v\in \RR^3$ 
in general position but close to the mirrors $H_{s_1}$ and $H_{s_3}$
the polytope $P(W;V,v)$ is pictured in Figure \ref{ExampleWythoff}.a).
For $V=\{s_2, s_3\}$ 
and $v\in H_{s_1}$ the polytope $P(W;V,v)$ is pictured
in Figure \ref{ExampleWythoff}.b).

\begin{proposition} \label{cox}
The combinatorial type of $P(W;V,v)$ is independent of the choice of $v$.
\end{proposition}
\proof For $V=S$ the polytope obtained is the well-known permutahedron, whose
face-lattice is independent of $v$.
The stabilizer of a face of $P(W; V, v)$ is a parabolic subgroup and
this establishes an isomorphism between the face lattice of $P(W; V, v)$
and the lattice of parabolic subgroups of $W$.
Furthermore, the $1$-skeleton of $P(W;S,v)$
is the Cayley graph $\Cay(W, S)$ of $W$ with respect to the generating set $S$.
Observe that in $\Cay(W, S)$ the length of any
edge labelled by generator $s\in S$
decreases as $v$ is moved towards the hyperplane $H_s$, and that
the edge ceases to exist when $v$ moves into $H_s$.
Denote by $\Cay_V(W, S)$ the obtained reduced Cayley graph.
The group $W$ acts transitively on its vertex set.
The set of vectors $v$ of type $t(v)=S$ has measure one in the set of all
vectors. Thus any face $F$ of $P(W;V,v)$ is obtained as the limit of some face
$G$ of the permutahedron and the set of vertices contained in $F$ is 
the limit of the vertices contained in $G$. Thus the graph $\Cay_V(W, S)$
determines the face-lattice and this proves the required result. \qed

We are going to describe explicitly the face lattice $P_{W, V}$
of the polytope $P(W;V,v)$ for any $v$ of type $V$.
Our description
is of course equivalent to the classic one given in 
\cite{Cox35,Cox73,Wyt} and reproduces the one of \cite{Sch90}.
Since we assume that $W$ is finite, the Coxeter-Dynkin reduced diagram
of $W$ is a tree and given any two vertices $u$ and $v$ of it
we denote by $[u,v]$ the unique path from $u$ to $v$.

For two subsets $U,U'\in S$ we say that $U'$
{\em blocks} $U$ (from $V$) if for all $u\in U$ and $v\in V$ there is
a $u'\in U'$, such that $u'\in [u,v]$.
This defines a binary relation on subsets of $S$, which we will
denote by $U'\le U$.
We also write $U'\sim U$ if $U'\le U$ and $U\le U'$, and we write $U'<U$ if
$U'\le U$ and $U\not\le U'$.

It is easy to see that $\le$ is reflexive and transitive, which implies
that $\sim$ is an equivalence relation. Let $[U]$ denote the equivalence
class containing $U$.
It can be shown that if $U\sim U'$
then $U\cap U'\sim U\sim U\cup U'$. This yields that every equivalence
class $X$ contains a unique smallest (under inclusion) subset $m(X)$ and
unique largest subset $M(X)$.
The subsets $m(X)$ will be called the {\em essential} subsets of $S$ (with
respect to $V$). Let $E(V)$ be the set of all essential subsets of
$S$. Clearly, the above relation $<$ is a partial order on $E(V)$. Also,
$V\in E(V)$ and $V$ is the smallest element of $E(V)$ with respect to $<$.

The faces $F$ of $P(W;V, v)$ are indexed by their isobarycenters $g(F)$.
The stabilizer of $F$ is the stabilizer of $g(F)$, that is the parabolic
subgroup of $W$ generated by $S - t(g(F))$.
The type of such an isobarycenter is an essential subset of $S$
and all essential subsets are realized as isobarycenters of faces.
The rank of
an essential subset is the dimension of the corresponding face.
Given two faces $F, F'$ of $P(W;V,v)$, $F\subset F'$ if and only if
we have the type inequality $t(g(F)) < t(g(F'))$ and
$\{g(F), g(F')\}$ is contained in at least one image $g(\Delta)$ with
$g\in W$ of the fundamental simplex $\Delta$.

We can use the above formalism to
obtain the combinatorial structure of the orbit polytope
$P(M_{24},v)$ where the Mathieu group acts on $\RR^{24}$ by permuting basis vectors, and $v=(1,2,3,4,5,0,\ldots,0)\in \RR^{24}$.
Since $M_{24}$ is a $5$-transitive permutation group we have
\begin{equation*}
P(M_{24},v) = P(S_{24},v).
\end{equation*}
The symmetric group $S_{24}$ is a finite reflection group with simple generating system $S=\{s_i=(i,i+1) \colon 1\le i\le 23\}$.
The vector $v$ lies in those mirrors $H_{s_i}$ for $6\le i\le 23$.
So $P(M_{24},v) = P(S_{24},V, v)$ for $V=\{s_1, \dots, s_5\}$.

Our proof of Proposition \ref{cox} implies that the polytope $P(M_{24},v)$
has $|S_{24}/\langle s_i \colon 6\le i\le 23\rangle| = 5100480$ vertices.
The essential subsets of rank $1$ defining edges are $V-\{s_k\}$ for $1\leq k\leq 4$ and $(V-\{s_5\})\cup \{s_6\}$. So, the number of edges is
\begin{equation*}
|S_{24}/\langle s_5, s_i:7\le i\le 23\rangle| +
\sum_{1\le k\le 4}|S_{24}/\langle s_k, s_i:6\le i\le 23\rangle| = 58655520.
\end{equation*}
Each vertex of the polytope has the same degree $d$ say.
Thus the number of edges is 
$d \times 5100480/2 = 58655520$ from which $d=23$.
Since $P(M_{24},v)$ is of dimension $23$, this shows that it is simple.

Each vertex of $P(M_{24},v)$ has stabilizer group $\Stab(M_{24},v)=M_{24}\cap \langle s_i:6\le i\le 23\rangle \cong 
(C_2 \times C_2 \times C_2 \times C_2) : C_3$ of order $48$.
Under $M_{24}$, for $1\leq k\leq 4$, there is only one orbit of edges of type $V-\{s_k\}$;
they have stabilizer $\Stab(M_{24},v) : C_2$ of order $96$.
Under $M_{24}$ there are two orbits of edges of type $(V-\{s_5\})\cup \{s_6\}$,
one with stabilizer $S_3$, the other with stabilizer a $2$-group of order $32$.

\medskip
The formalism of essential subsets is a useful tool to determine
the face lattice of $P (W ; V , v )$ for a Coxeter group $W$
and provides ready access to the lattice for homology computations.
The equality between the polytopes $P(M_{24},v)$ and $P(S_{24}, v)$ 
was essential for being able to apply this formalism and thus get
a reasonably simple description of the face lattice.

For an arbitrary vector $v$ and group $G$ we cannot expect
to have a simple combinatorial description of the face lattice of $P(G, v)$
and we
need to use specific computational techniques.
If $G$ is large, then we cannot
expect to be able to store the vertex set of $P(G, v)$.
Fortunately, by the group action, the full face lattice is
encoded in the set $S(v)$ of vertices adjacent to $v$.
This set $S(v)$ can be computed iteratively by using the
Poincar\'e polyhedron theorem (see \cite{Riley,deraux} for some example
of such computations).
Once the list of neighbours is known the face-lattice follows easily.

After one has obtained the low dimensional faces of $P(M_{24}, v)$ and
their stabilizer groups, we can use Theorem \ref{wall} to compute the
initial terms of a free $\ZZ M_{24}$-resolution of $\ZZ$.

\section{Wythoff construction for polytopes}

The Wythoff construction can also be defined for partially 
ordered sets.
A {\em flag} in a poset is an arbitrary completely ordered subset.
We say that a connected poset $\cK$ is a {\em $d$-dimensional complex} (or,
simply, a {\em $d$-complex}) if every maximal flag in $\cK$ has size
$d+1$. In a $d$-complex $\cK$ every element $x$ can be uniquely assigned a
number $\dim(x)\in\{0,\ldots,d\}$, called the {\em dimension} of $x$, in
such a way, that the minimal elements of $\cK$ have dimension zero and
$\dim(y)=\dim(x)+1$ whenever $x<y$ and there is no $z$ with $x<z<y$.
The elements of a complex $\cK$ are called {\em faces}, or {\em $k$-faces}
if the dimension of the face needs to be specified.  Furthermore, $0$-faces
are called {\em vertices} and $d$-faces (maximal faces) are called {\em
facets}. If $x<y$ and
$\dim(x)=k$, we will say that {\em $x$ is a $k$-face of $y$}.

For a flag $f\subset\cK$ define its {\em type} as the set
$t(f)=\{\dim(F)\quad : \quad F\in f\}$. Clearly, $t(f)$ is a subset of
$S=\{0,\ldots,d\}$ and, conversely, every subset of $S$ is the type of
some flag.
Let $\Omega$ be the set of all nonempty subsets of $S$ and fix an
arbitrary $V\in\Omega$. For two subsets $U,U'\in\Omega$ we say that $U'$
{\em blocks} $U$ (from $V$) if for all $u\in U$ and $v\in V$ there is
a $u'\in U'$ and $u\leq u'\leq v$ or $v\leq u'\leq u$.
With this notion of blocking we can define the notion of essential subset
of $S$ and the inequality $<$ in the same way as for Coxeter groups.

The construction of $P(\cK;V)$ mimics the one of $P(W;D,v)$ above 
for Coxeter groups.
The {\em Wythoff complex} $P(\cK;V)$
consists of all flags $F$ such that $t(F)$ is essential.
For two such flags $F$ and $F'$, we have $F'<F$ whenever $t(F')<t(F)$ and $F'$
is compatible with $F$, that is, $F\cup F'$ is a flag.
It can be shown that $P(\cK, V)$ is again a $d$-complex.

The face lattice $\cK(P)$ of a $(d+1)$-dimensional polytope $P$ is a
$d$-complex, which is a CW-complex
topologically equivalent to a sphere.
It is proved in \cite{Sch90} that the topological type of $P(\cK;V)$
is the same as the one of $\cK$.
This version of the Wythoff construction when applied to a regular polytope 
gives a face lattice which is isomorphic to the one obtained by applying
the Wythoff construction to the corresponding Coxeter group.
The complex $P(\cK(P), \{0\})$ is equal to $\cK(P)$ and $P(\cK(P), \{d\})$
is the complex of the polytope dual to $P$.
In general $P(\cK(P), V)$ is not a polytope since the notion
of convexity is not well preserved by the Wythoff construction
without any regularity assumption.

The topological invariance means that if a group $G$ acts on a polytope
$P$ then we can apply the orbit polytope construction to $P(\cK(P), V)$
for a chosen $V$ in order to compute $H_i(G, \ZZ)$.

In the case of $M_{24}$, we take as polytope the $23$-dimensional simplex
$\alpha_{23}$ and we build the Wythoff polytope
$P(\alpha_{23}; \{0,1,2,3,4\})$.
In Table \ref{WythoffResolutions} we give the results obtained for the
larger Mathieu groups.
The method applies to any finite group acting on $n$ points by using the
simplex $\alpha_{n-1}$.
We do not need $G$ to act transitively.
All programs are available from \cite{polyhedral}.

\begin{table}
\begin{center}
\begin{tabular}{c|c|c|l}
$G$      & $P$          & $V$ & Free rank of resolution in degrees $0,1,2,\ldots$\\
\hline
$M_{22}$ & $\alpha_{21}$ & $\{0,1,2\}$ & $1$, $7$, $33$, $113$, $301$, $694$\\
$M_{23}$ & $\alpha_{22}$ & $\{0,1,2,3,4\}$ & $2$, $20$, $116$, $451$, $1334$, $3279$\\
$M_{24}$ & $\alpha_{23}$ & $\{0,1,2,3,4\}$ & $1$, $9$, $50$, $204$, $649$\\
\end{tabular}
\end{center}
\caption{Rank of resolutions of $M_{22}$, $M_{23}$, $M_{24}$ obtained from the Wythoff construction}
\label{WythoffResolutions}
\end{table}

\section{Homology at $p=5,7,11,23$}
Suppose that a group $G$ has Sylow $p$-subgroup $P=C_p$ of prime order.
The Cartan-Eilenberg double coset formula implies that the surjection 
$$\pi_n\colon H_n(P,\ZZ)\rightarrow H_n(G,\ZZ)_{(p)}$$
has kernel generated by the elements
$$ H_n(\phi_g)(a) -a $$
for $ g\in N_G(P), a\in H_n(P,\ZZ)$ and $\phi_g\colon P
\rightarrow P,p\mapsto gpg^{-1}$.
Here $N_G(P)$ is the normalizer of $P$ in $G$.

Using the isomorphism $H_{n-1}(P,\ZZ)\cong H^n(P,\ZZ)$ and 
the cohomology ring structure $H^\ast(P,\ZZ) \cong \ZZ_p[x^2]$,
we see that a group homomorphism $\phi\colon P\rightarrow P, p\mapsto p^m$
induces a homology homomorphism $H_{2k-1}(\phi)\colon H_{2k-1}(P,\ZZ) 
\rightarrow H_{2k-1}(P,\ZZ), a\mapsto a^{m^k}$.

For $p\in \{5,7,11,23\}$ the Mathieu groups have Sylow $p$-subgroups which are either trivial or of prime order.
One can use {\sc gap} to determine their normalizers.
It is thus a routine exercise to determine the
$p$-part of the integral homology of the Mathieu groups, the 
results of which are given in the Introduction.

\end{document}